\documentclass[11pt,a4paper]{article}
\usepackage{epsfig}
\usepackage[T1]{fontenc}    % Accents cods dans la fonte.
\usepackage{graphics}
\usepackage{graphicx}
\usepackage{pstricks,pst-coil,pst-fill,pst-plot}
\usepackage[fleqn]{amsmath}    % Les symboles les plus frquents.
\usepackage{amssymb}    % Des symboles.
\usepackage{amsfonts}   % Des fontes, eg pour \mathbb.
\usepackage{verbatim}   % Pour les codes sources en informatique.
\usepackage{mathrsfs}   % Des lettres majuscules cursives (\mathscr).
\usepackage{dsfont}
\usepackage{euscript}
\usepackage{yfonts}
\usepackage{enumerate}     % met en place l'environement enumerate pour les listes
\usepackage{txfonts}
\usepackage{marvosym}
\usepackage{vmargin}        % Rgler la taille de la feuille.

\setmarginsrb{1.8cm}{2cm}{1.8cm}{2cm}{1cm}{1cm}{1cm}{1.6cm}
 \makeatletter
 \@addtoreset{equation}{section}
 \makeatother

\providecommand{\bysame}{\leavevmode\hbox to3em{\hrulefill}\thinspace}
\providecommand{\MR}{\relax\ifhmode\unskip\space\fi MR }
\providecommand{\href}[2]{#2}

\let\tend=\rightarrow

\long\def\symbolfootnote[#1]#2{\begingroup%
\def\thefootnote{\fnsymbol{footnote}}\footnote[#1]{#2}\endgroup}

\newtheorem{theorem}{Theorem}[section]
\newtheorem{prop}{Proposition}[section]

\def\Proof{\medskip\noindent {\it Proof --- \ }}

\def\qed{\hfill\rule{2mm}{2mm}}

\newcommand\beq{\begin{equation}}
\newcommand\enq{\end{equation}}
\newcommand\bem{\begin{multline}}
\newcommand\enm{\end{multline}}

\def\beqa{\begin{eqnarray}}
\def\eeqa{\end{eqnarray}}
\def\ba{\begin{array}}
\def\ea{\end{array}}
\def\det{\operatorname{det}}

\newcommand{\f}[2]{{\ensuremath{%
    \mathchoice%
    {\dfrac{#1}{#2}}
    {\dfrac{#1}{#2}}
    {\frac{#1}{#2}}
    {\frac{#1}{#2}}
}}}
\newcommand{\tf}[2]{\ensuremath{#1/#2}}
\newcommand{\pa}[1]{\ensuremath{\left(#1\right)}}
\def\a{\alpha}

\def\ga{\gamma}
\def\Ga{\Gamma}

\def\veps{\varepsilon}
\def\la{\lambda}

\def\sg{\sigma}

\newcommand{\mc}[1]{\ensuremath{\mathcal{#1}}}

\newcommand{\msc}[1]{\ensuremath{\mathscr{#1}}}

\newcommand{\bs}[1]{\ensuremath{\boldsymbol{#1}}}
\newcommand{\ov}[1]{\ensuremath{\overline{#1}}}
\newcommand{\wt}[1]{\ensuremath{\widetilde{#1}}}

\newcommand{\Int}[2]{\ensuremath{\int\limits_{#1}^{#2}}}
\newcommand{\Oint}[2]{\ensuremath{\oint\limits_{#1}^{#2}}}

\newcommand{\sul}[2]{\ensuremath{\sum\limits_{#1}^{#2}}}
\newcommand{\pl}[2]{\ensuremath{\prod\limits_{#1}^{#2}}}

\newcommand{\R}{\ensuremath{\mathbb{R}}}
\newcommand{\Cx}{\ensuremath{\mathbb{C}}}

\newcommand{\Dp}[1]{\ensuremath{\partial_{#1}}}
\newcommand{\ex}[1]{\ensuremath{\e{e}^{#1}}}

\newcommand{\abs}[1]{\ensuremath{\left| #1 \right|}}

\newcommand{\norm}[1]{\ensuremath{|| #1 ||}}

\newcommand{\dd}{\mathrm{d}}
\newcommand{\e}[1]{\ensuremath{\mathrm{#1}}}

\newcommand{\intff}[2]{\ensuremath{ [  #1 \,; #2 ] }}

\newcommand{\intn}[2]{\ensuremath{[\![ \, #1 \,;\, #2 \,]\!]}}

\begin{document}

\begin{flushright}

\end{flushright}
\par \vskip .1in \noindent

\vspace{14pt}

\begin{center}
\begin{LARGE}
{\bf On determinants of integrable operators with shifts.}
\end{LARGE}

\vspace{30pt}

\begin{large}

{\bf A.~R.~Its}\symbolfootnote[1]{Indiana Univercity Purdue Univercity Indianapolis, Department of Mathematics, Indianapolis, USA, itsa@math.iupui.edu} , 
{\bf K.~K.~Kozlowski}\symbolfootnote[2]{Universit\'{e} de Bourgogne, Institut de Math\'{e}matiques de Bourgogne, UMR 5584 du CNRS, France,
karol.kozlowski@u-bourgogne.fr}. 
\par

\end{large}

\vspace{40pt}

\centerline{\bf Abstract} \vspace{1cm}
\parbox{12cm}{\small 
Integrable integral operator can be studied by means of a matrix Riemann--Hilbert problem. 
However, in the case of so-called integrable operators with shifts, the associated Riemann--Hilbert problem 
becomes operator valued and this complicates strongly the analysis. 
In this note, we show how to circumvent, in a very simple way,
the use of such a setting while still being able to characterize the large-$x$ asymptotic behavior of the 
determinant associated with the operator.}

\end{center}

\vspace{40pt}

\section*{Introduction}

The theory of integrable integral operators takes its roots in the works of Jimbo, Miwa, Mori and Sato 
\cite{JimMiwaMoriSatoSineKernelPVForBoseGaz} which ultimately led to the development\footnote{
Several principal aspects of the integrable operator theory, especially
the ones concerning
 the integrable differential systems appearing in random matrix theory,
have been developed in \cite{TracyWidomPDEforFredholms}. Some of the  elements of
the  theory of integrable operators were already implicitly present in the
earlier work \cite{SakhnovichFirstStepsOfElementsOfIntegrableIntOps}.}
\cite{ItsIzerginKorepinSlavnovDifferentialeqnsforCorrelationfunctions} of a
Riemann--Hilbert based setting for studying these operators. The initial motivation
for studying these operators stemmed from the theory of quantum integrable models at the free fermion point
where the correlation functions are expressed in terms of Fredholm determinants -or minors thereof-. 
Let us recall certain basic facts. An integrable integral operator is an integral operator on $L^2\big(\msc{C}\big)$
of the type $I+\wt{V}$ where $\msc{C}$ is some contour in $\Cx$ and the integral kernel $\wt{V}$ takes the form\footnote{the last condition may be avoided for the price of dealing with a principal-value regularization, see \cite{DeiftIntegrableOperatorsDiscussion}}

\beq
\wt{V}(\la,\mu) \; = \; \f{ \sum_{a=1}^{N} f_a(\la) e_a(\mu) }{ \la - \mu }  \qquad 
\e{with} \qquad \sul{a=1}{N}  f_a(\la) e_a(\la) \; = \; 0 \;. 
 \label{definition integrable kernel}
\enq
There is no restriction on the number $N$ of functions involved in the definition of the kernel. 
It has been shown in \cite{ItsIzerginKorepinSlavnovDifferentialeqnsforCorrelationfunctions} that the study of such integral 
operators (calculation of their resolvent kernel, determinant, \textit{etc} ) boils down 
to the resolution of an $N\times N$ Riemann--Hilbert problem. This approach is particularly efficient
in the asymptotic regime when the functions $f_a$ and $e_a$ entering in the composition of the kernel depend on 
some large parameter $x$, say in an oscillatory way. Indeed, then, the jump matrix arising in the associated 
Riemann--Hilbert problem depends on the large-parameter $x$ in an oscillatory way and its asymptotic solution can be studied within the Deift-Zhou non-linear steepest descent method \cite{DeiftZhouSteepestDescentForOscillatoryRHP}. This program has been carried out 
for a large number of kernels in the $N=2$ case, see \textit{eg} \cite{DeiftItsZhouSineKernelOnUnionOfIntervals,KozKitMailSlaTerRHPapproachtoSuperSineKernel}. 

In fact, one can even consider more general kernels where the discreet sum in \eqref{definition integrable kernel}
is replaced by an integration in respect to an arbitrary measure $\dd \mu$ on some set $E$:
\beq
V(\la,\mu) \; = \; \f{ \int_{E} f(\la,s) \cdot e(\mu,s)  \cdot \dd \mu(s)  }{ \la - \mu }  \qquad 
\e{under} \; \e{the} \; \e{condition} \qquad 
\int_{E}  f(\la,s) \cdot e(\la,s) \cdot \dd \mu(s)  \; = \; 0 \;. 
 \label{definition integrable kernel general}
\enq
There is however a price to pay for such a generalization: as soon as the measure $\dd \mu$ is not 
compactly supported and purely atomic, one deals with an operator valued Riemann--Hilbert problem. 
In the early days of the analysis of integrable integral operators, there have been a few
attempts to extract valuable data out of such Riemann--Hilbert problem. 
The matter is that, even in the most simple cases, the analysis becomes extremely hard and complicated due to the operator valued setting
for the Riemann--Hilbert problem. 
Still, the authors of \cite{ItsSlavnovNLSTimeAndSpaceCorrDualFields} were able to 
argue the leading asymptotics of the solution to the Riemann--Hilbert problem arising in the large-$x$
asymptotic analysis of a very specific integrable integral operator of the form \eqref{definition integrable kernel general}
with a measure $\dd \mu(s)$ supported on $\R^+$ that was absolutely continuous in respect to Lebegue's one. 
The purpose of the present paper is to develop
an alternative approach to the study of operators of the type \eqref{definition integrable kernel general}
for certain classes of functions $f$ and $e$ and in the context of specific measures $\dd \mu(s)$ supported on $\R^+$
which are absolutely continuous in respect to Lebesgue's one. 
For instance, the case considered in \cite{ItsSlavnovNLSTimeAndSpaceCorrDualFields} does fall into this class. 
As an application of our technique, we shall consider the operator $I+S$ acting on $L^2\big( \intff{a}{b}\big)$ with an 
integrable integral kernel $S(\la,\mu)$ of, so-called, shift type: 
\beq
S(\la,\mu) \; = \; \f{ ic F(\la) }{ 2i\pi (\la - \mu) } 
\bigg\{  \f{ e(\la) e^{-1}(\mu)  }
{ (\la-\mu + i c) }  \; + \; 
\f{ e(\mu) e^{-1}(\la) }
{ (\la - \mu - i c) } \bigg\}  \qquad \e{where} \qquad e(\la) \; = \; \ex{ i \f{x}{2}p(\la) } \;.
\label{definition noyau type shift}
\enq
$F$ and $p$ are certain holomorphic function in a neighborhood of the interval $\intff{a}{b}$
and $F$ is "sufficiently" small, in a sense that will be specified below.  
We call these kernels of "shift-type" in that the denominator not only contains the singular factor $\la-\mu$
but also shifts thereof,  $\la-\mu \pm ic$ in this case. 
It is readily seen that the kernel \eqref{definition noyau type shift}
is an integrable kernel in the generalized sense \eqref{definition integrable kernel general} since it 
admits the representation
\beq
S(\la,\mu) \; = \;  \f{-ic}{(\la-\mu)}  \Int{0}{+\infty} \f{F(\la) \ex{-c s} }{2i\pi} 
\big[ e(\la) \cdot e^{-1}(\mu)  \cdot \ex{is(\la-\mu)} \, - \,  e(\mu) \cdot e^{-1}(\la) \cdot  \ex{is(\mu-\la)}  \big] \cdot \dd s \;. 
\enq
Our approach allows us to compute the large-$x$ asymptotic behavior of the Fredholm determinant of the operator 
$I+S$, what constitutes the main result of the paper. 
\begin{theorem}
\label{Theorem principal}
Let $\wt{S}$ correspond to the below generalized sine kernel 
\beq
\wt{S}(\la,\mu)  \; = \;  F(\la) \f{  e(\la)e^{-1}(\mu) - e(\mu)e^{-1}(\la) }{ 2i\pi (\la- \mu) }   \;, 
\label{definition operateur S tilde noyau}
\enq
with $F$, $p$ holomorphic functions on some open neighborhood $U$ of $\intff{a}{b}$ and such that $|F|<1$ on $U$
and $p^{\prime}_{\mid \intff{a}{b}} >0$. Then, the Fredholm determinant of the integral operator $I+S$
acting on $L^2(\intff{a}{b})$ with a kernel $S(\la,\mu)$ given by \eqref{definition noyau type shift}
admits the $x \tend +\infty$ asymptotic expansion
\beq
\f{ \det{}\big[ I\, + \, S \big]  }{ \det{}\big[ I\, + \, \wt{S}\big]  }  \; = \; 
\underset{ \Ga }{\det}\big[ I\; + \; U_+  \big] \cdot \underset{ \Ga }{\det}\big[ I\; + \; U_-  \big]
\; \cdot \;\Big\{  1  \; + \; \e{O}\big( x^{-1} \big)   \Big\} \;. 
\enq
The asymptotic expansion of $\det{}\big[ I\, + \, \wt{S}\big]$ can be found in \cite{KozKitMailSlaTerRHPapproachtoSuperSineKernel}
and $I+U_{\pm}$ are integral operators acting on $L^2\big( \Ga \big)$  with kernels 
\beq
U_{+}\big(\la, \mu \big) \; = \;    \f{ \a(\mu - ic)\cdot \a^{ -1 }(\la) }{ 2i\pi \big(\la - \mu + ic \big) }
\qquad \e{and} \qquad U_{ - }\big(\la, \mu \big) \; = \;   \f{ \a^{- 1}(\mu + ic) \cdot \a(\la) }{ 2i\pi \big(\la - \mu - i c \big) } \;, 
\enq
whereas $\Ga$ is a small counterclockwise loop around $\intff{a}{b}$ of index 1. 
Finally, the function $\a$ appearing above is given by 
\beq
\a(\la) \; = \; \exp\Bigg\{ \Int{a}{b}  \f{\ln [1+F(\mu)] }{ \la- \mu } \cdot  \f{ \dd \mu }{ 2i\pi } \Bigg\} \;. 
\label{definition fonction alpha}
\enq

\end{theorem}

This theorem can be as well obtained by the method of multidimensional Natte series whose premises appeared in 
\cite{KozKitMailSlaTerXXZsgZsgZAsymptotics}
and which has been further developped and simplified in \cite{KozReducedDensityMatrixAsymptNLSE,KozTerNatteSeriesNLSECurrentCurrent}. 
In particular, the method has been brought to a satisfactory level of rigor in 
\cite{KozReducedDensityMatrixAsymptNLSE} upon a hypothesis of convergence of a series of multiple integrals. 
In the case underlying to the kernel $S$, 
it is fairly easy to establish this convergence. However, the multidimensional Natte series method is definitely much more complicated
that the one being developed in the present note, so that 
applying it to the analysis of such "simple operators" would be, mildly speaking, highly unreasonable.

We also stress that our method is quite general and not solely applicable to the above kernel. As will be apparent from the core of the text, 
it allows one even to obtain new types of representations for Fredholm determinants of non-integrable operators of the type
\beq
V(\la,\mu) \; = \;   \wt{V}(\la,\mu)  \; - \; 
\sul{a=1}{N} \ga_a \f{ f_a(\la) e_{v_a}(\mu) }{ \la-\mu + i c_a } \qquad 
\e{where}  \qquad \wt{V}(\la,\mu) \; = \; \f{ \sum_{a=1}^{N} f_a(\la) e_a(\mu) }{ \la - \mu } \; , 
\label{definition noyau avec shifts}
\enq
$c_a \in  \R^{*}$ and $v_a$ is an arbitrary sequence in $\intn{1}{N}$.  For generic constants $\ga_a$ and sequences $v_a$,  $a=1,\dots,N$
the operators' kernels are not of integrable type, but they reduce to the integrable case for the specific choice $\ga_a=1$
and $v_a = a$ for any $a$. In fact, in such an integrable integral operator case, it would be interesting to compare 
the approach developed in the present note with the results that could be obtained with the 
help of operator valued Riemann-Hilbert techniques of \cite{ItsSlavnovNLSTimeAndSpaceCorrDualFields}. 
We are going to address this issue in a forthcoming publication.

The paper is organized as follows. In section \ref{Section RHP setting}, we recall several basic facts about Riemann--Hilbert problems for 
integrable integral operators. Then, in section \ref{Section Factorization of determinants}, we establish a factorization for the determinant $\det[I+V]$ in terms of $\det[I+\wt{V}]$ and antother determinant involving the solution to the 
Riemann--Hilbert problem associated with the integrable kernel $\wt{V}$. Finally, in section \ref{Section Application of Results} 
we demonstrate that in the case of the integrable operator of shift-type \eqref{definition noyau type shift}
associated with the generalized sine kernel, such a factorization is already enough 
so as to access to the asymptotics of the determinant $\det[I+S]$ which take the form stated in theorem \ref{Theorem principal}.

\section{The Riemann--Hilbert problem setting }
\label{Section RHP setting}

\subsection{The kernel and initial Riemann--Hilbert problem}

Let $I+ \wt{V}$ be an integral operator acting on $L^{2}\big( J \big)$ with $J$ a piecewise smooth curve in $\Cx$
and  whose kernel $\wt{V}(\la,\mu)$ 
is of integrable integral type :
\beq
\wt{V}(\la,\mu) =  \f{ \big( \bs{E}_{L}(\la) , \bs{E}_{R}(\mu) \big) }{ \la-\mu}   
\enq
with vector valued smooth functions on $J$
\beq
 \bs{E}_{L}^{\bs{T}}(\la) \; = \; \big( f_1(\la), \dots, f_N(\la) \big) \qquad \e{and} \qquad  
\bs{E}_{R}^{\bs{T}}(\la) \; = \;  \big( e_1(\la), \dots, e_N(\la) \big) \; 
\enq
satisfying to the regularity condition $\big( \bs{E}_{L}(\la) , \bs{E}_{R}(\la) \big)=0$ for all $\la \in J$. 

This kernel is associated with the Riemann--Hilbert problem for a $N\times N$ matrix $\chi(\la)$
\begin{itemize}
\item $\chi \in \mc{O}(\Cx\setminus J )$ and has continuous boundary values on 
$\overset{ \circ }{ J }$;
\item $\chi(z)=\ln\abs{z-a}  \e{O}\pa{ \ba{cc} 1 & 1 \\ 1 & 1 \ea}$ 
as $z \tend a$ with $a \in  \Dp{} J $; 
\item $\chi(z) = I_{N} \; + \; \e{O}\big(z^{-1}\big)$ when $z \tend \infty$;
\item $\chi_-(z) = \chi_+(z) \cdot G_{\chi}(z)$ for $z \in \overset{ \circ}{ J }$ where 
$G_{\chi}(\la) = I_N + 2i\pi\,  \bs{E}_R^{\bs{T}}(\la)\cdot \bs{E}_{L}(\la)$  \;. 
\end{itemize}
Here, we should explain that relations of the type $M(z) = \e{O}(R(z))$ for two matrix functions $M,R$ should be understood entry-wise, 
\textit{ie} $M_{jk}(z)= \e{O}\big( R_{jk}(z) \big)$. Also, given a function $f$ defined on $\Cx\setminus \ga$, with $\gamma$
an oriented curve in $\Cx$, we denote by $f_{\pm}(s)$  the boundary values of $f(z)$ on $\ga$ when the argument $z$ 
 approaches the point $s \in \ga$ non-tangentially and from the left ($+$) or the right ($-$) side of the curve. Again, if one deals with 
matrix function, then this relation has to be understood entry-wise. 
 
The unique solvability of the above Riemann--Hilbert problem for $\chi$ is equivalent to the condition $\det[I+\wt{V}] \not=0$. 
In such a case, this unique solution takes the form
\beq
\chi(\la) = I_N - \Int{J}{}  \f{ \bs{F}_R(\mu) \cdot \bs{E}_L^{\bs{T}}(\mu) }{ \mu- \la }  \dd \mu    \qquad \e{and} \qquad
\chi^{-1}(\la) = I_N + \Int{J}{}  \f{ \bs{E}_R(\mu) \cdot \bs{F}_L^{\bs{T}}(\mu)  }{ \mu - \la }  \dd \mu 
\label{formules reconstruction chi chi-1 en terms F R et FL}
\enq
where $I_N$ is the $N\times N$ identity matrix and $\bs{F}_{R}(\mu)$ and $\bs{F}_{L}(\mu)$ correspond to the solutions to the below linear integral equations
\beq
\bs{F}_{R}(\la)  \; + \; \Int{J}{} \wt{V}(\mu,\la) \bs{F}_{R}(\mu)  \cdot  \dd \mu \; = \;   \bs{E}_R(\la)
\qquad \e{and} \qquad 
\bs{F}_L(\la)  \; + \; \Int{J}{} \wt{V}(\la,\mu) \bs{F}_L(\mu) \cdot  \dd \mu\; = \;   \bs{E}_L(\la) \;. 
\enq

Also, the vector functions  $\bs{F}_{R}(\mu)$ and $\bs{F}_{L}(\mu)$ can be reconstructed in terms of $\chi$
as follows
\beq
\bs{F}_{R}(\mu) \; = \; \chi(\la) \cdot  \bs{E}_{R}(\mu)\qquad \e{and}  \qquad \bs{F}_{L}^{\bs{T}}(\mu) \; = \; 
 \bs{E}_{L}^{\bs{T}}(\mu) \cdot \chi^{-1}(\la) \;. 
\enq

\section{Factorization of determinants}
\label{Section Factorization of determinants}

In this section, we consider a trace class integral operator $I+V$ on $L^2(J)$ with a kernel given by 
\eqref{definition noyau avec shifts}. We establish various alternative representations for its 
Fredholm determinant. 

\begin{prop}
\label{Proposition factorization determiant V en terme RHP sovable simple det V tile et det W}

Assume that $ \det[I+ \wt{V} ] \not= 0$. Then, the Fredholm determinant $\det[I+V]$ admits the decomposition
\beq
\det[I+V] \; = \; \det[I+ \wt{V} ] \cdot \det[I+W]
\enq
where the integral kernel $W(\la,\mu)$ takes the form 
\beq
W(\la,\mu) \; = \; - \sul{ n = 1 }{ N } \ga_n \f{ \big(\bs{F}_{L}(\la) ,  \chi(\mu \,- i c_n ) \,  \bs{e}_n \big) 
												\cdot \big( \bs{e}_{v_n}, \bs{E}_{R}(\mu) \big) }
			{ \la \; - \; \mu  \, + \,  i c_n   }			\; , 
\enq
in which we have introduced  
\beq
\bs{e}_p^{\bs{T}}   \; = \; \big( \underbrace{0, \dots, 0}_{p-1}, 1, 0,\dots , 0 \big) \; \in \R^N \;.
\enq

\end{prop}

\Proof 

It follows that 
\beq
I+V \; = \;  \big( I + \wt{V} \big) \cdot \big(I + W \big)
\enq
where the kernel 
\beq
W(\la,\mu)  \; = \; -  \sul{a=1}{N} \ga_a  \mc{U}_a(\la,\mu) e_{v_a}(\mu)  
\enq
is expressed in terms of the unique solutions to the integral equations
\beq
\mc{U}_n(\la,\mu) \; + \; \Int{J }{  } \wt{V}(\la,\tau) \mc{U}_n(\tau,\mu) \cdot \dd \tau  \; = \; 
\f{ \big( \bs{E}_L(\la), \bs{e}_n \big) }{ \la - \mu + i c_n } \;. 
\enq
One can solve these equations in terms of the entries of $\chi$, the unique solution of the RHP associated with 
the kernel $\wt{V}(\la,\mu)$. Setting $\msc{U}_n(\la,\mu)  \, = \,  (\la-\mu +  i c_n) \mc{U}_n(\la,\mu)$, 
we get that $\msc{U}_{n}(\la,\mu)$
solves the integral equation 
\beq
\msc{U}_{n}(\la,\mu) \; + \; \Int{ J }{ } \wt{V}(\la,\tau) \msc{U}_{n}(\tau,\mu)  \cdot  \dd \tau  
 \; + \; \Int{ J }{ } \big( \bs{E}_L(\la) , \bs{E}_{R}(\tau) \big) \cdot \mc{U}_{n}(\tau,\mu) \cdot \dd \tau \; = \; 
 \big( \bs{E}_L(\la), \bs{e}_n \big)   \;. 
\enq
Thus one has that
\beq
\mc{U}_{n}(\la,\mu) \; = \; \f{  \big( \bs{F}_L(\la), \bs{e}_n \big)  }{ \la \, - \,  \mu \, + \,  i c_n  }
\; - \;  \Int{ J }{ }  \f{\big( \bs{F}_L(\la) , \bs{E}_{R}(\tau) \big)  }{ \la-\mu \, +  \, ic_n  }
\cdot \mc{U}_{n}(\tau,\mu) \cdot \dd \tau  \;. 
\enq
It solely remains to fix the unknown vector coefficient. This can be done by solving the consistency conditions
\beq
\Int{ J }{ } \bs{E}_{R}(\tau) \cdot \mc{U}_{n}(\tau,\mu) \dd \tau \; = \; 
\Bigg( -  \Int{ J }{ }  \f{  \bs{E}_{R}(\tau) \cdot \bs{F}_{L}^{\bs{T}}(\tau) }{ \tau-\mu \, + \, i c_n }  \dd \tau \Bigg) 
\; \cdot \;  \Bigg( \Int{ J }{ } \bs{E}_{R}(\tau) \cdot \mc{U}_n(\tau,\mu) \dd \tau  \Bigg)
\;+\;   \Int{ J }{ }  \f{ \bs{E}_{R}(\tau) \cdot \bs{F}_{L}^{\bs{T}}(\tau) \cdot \bs{e}_n  }{ \tau-\mu \, + \, i c_n }   \;. 
\enq
 Using the integral representation for $\chi^{-1}(\la)$, we are let to the equation
\beq
\chi^{-1}(\mu \, -  i c_n ) \Bigg( \Int{ J }{ } \bs{E}_{R}(\tau) \cdot \mc{U}_n(\tau,\mu) \cdot  \dd \tau  \Bigg) \; = \; 
\Big[ \chi^{-1}(\mu \, - i c_n ) - I_N \Big] \cdot  \bs {e}_n
\enq
 which is readily solved for the unknown vector coefficients. Therefore
\beq
\mc{U}_n(\la,\mu) \; = \; 
 \f{ \big( \bs{F}_{L}(\la) ,   \chi(\mu \, - i c_n )  \bs{e}_n \big)  }{ \la \; - \; \mu  \, + \, i c_n   } \;. 
\enq
Upon  inserting these expression into $W(\la,\mu)$, we are led to the claim. \qed

The determinant of $I+W$ can be recast in terms of the Fredholm determinant of a matrix integral operator
acting on a small counterclockwise loop  $\Ga(J)$ encircling the original interval $J$.

\begin{prop}
\label{Proposition factorisation Det W en terme RHP data det M}
One has the representation
\beq
\underset{J}{\det} \big[I+W \big] \; = \; \underset{ \Ga(J) }{\det} \big[ I+ \mc{M} \big]
\enq
Where $\mc{M}$ is a $N\times N$ matrix integral operator on $L^2\big(  \Ga(J) \big)\otimes \Cx^N$ whose entries are given by 
\beq
\mc{M}_{k\ell}(\la,\mu) \; = \;    
 \ga_k \f{ \big( \bs{e}_{v_k} , \chi^{-1}(\la) \cdot \chi(\mu - i c_{\ell} )  \cdot \bs{e}_{\ell} \big) }
 { 2i\pi ( \la-\mu +i c_{\ell} ) }  \;. 
\label{ecriture matrice M}
\enq
\end{prop}

\Proof 

Let $\kappa \in \Cx$ and $|\kappa|$ be small enough, then one has the series expansion  
\beq
\ln \underset{J}{\det} [I+\kappa W] \; = \; \sul{ n \geq 1}{} \f{ (-1)^{n-1} }{ n } (-\kappa)^n  
 \sul{ \substack{ \tau_p \in \mc{E}_N  \\ p=1,\dots, n  }  }{} \pl{a=1}{n} \{ \ga_{\tau_a} \} \cdot  I_{\{\tau_p\}}^{(n)}
\qquad \e{where} \qquad \mc{E}_N \; = \; \intn{1}{N}
\enq
and 
\beq
I_{\{\tau_p\}}^{(n)} \; = \; \Int{ J }{ }  \hspace{-1mm}
\pl{p=1}{n}  \Bigg\{ \f{ \big( \bs{F}_{L}(\la_p),  \, \chi(\la_{p+1} \, - i  c_{\tau_p} ) \cdot \bs{e}_{\tau_p} \big) } 
{  \la_p \; - \; \la_{p+1}  \, + \,  i c_{ \tau_p }    }  
\cdot \big( \bs{e}_{v_{\tau_p}}, \bs{E}^R(\la_{p+1})  \big)  \Bigg\}\cdot  \dd^n \la \;. 
\enq
Above, we agree upon the cyclicity notation $\la_{n+1}=\la_1$. 
In order to be able to take the $\la$-integrals, we separate the $\la$-integration by introducing auxiliary $z$-contour integrals
which leads to 
\bem
I_{\{\tau_p\}}^{(n)} \; = \;
\Int{J}{} \dd^n \la  \Oint{\Ga(J)}{} \f{\dd^n z}{ (2i\pi)^n}
\pl{p=1}{n}  \f{  \big( \bs{F}_{L}(\la_p),  \, \chi(z_{p+1} \, - i  c_{\tau_p} ) \cdot \bs{e}_{\tau_p} \big) 
\cdot \big( \bs{e}_{v_{\tau_p}}, \bs{E}^R(\la_{p+1})  \big) } 
{  (z_p \; - \; z_{p+1}  \, + \,  i  c_{\tau_p} )(z_p-\la_p)   } \\
=   \Oint{\Ga(J)}{} \f{\dd^n z}{ (2i\pi)^{n}} 
\pl{p=1}{n} \f{ 1 }{ z_p \; - \; z_{p+1}  \, + \,  i c_{\tau_p} }
\Int{J}{} \dd^n \la  
\pl{p=1}{n}  
\f{ \big( \bs{e}_{v_{\tau_{p-1}}}, \bs{E}^R(\la_{p})  \big) \cdot \big( \bs{F}_{L}(\la_p),  \, \chi(z_{p+1} \, - i  c_{\tau_p} ) \cdot \bs{e}_{\tau_p} \big) } 
{ (z_p-\la_p)   }  \\
=  \Oint{\Ga(J)}{} 
\pl{p=1}{n}  \f{ \big( \bs{e}_{v_{\tau_{p-1}}},  \big[ I_N - \chi^{-1} (z_p) \big] \cdot 
\chi(z_{p+1} \, - i  c_{\tau_p} ) \cdot \bs{e}_{\tau_p} \big) } 
{ z_p \; - \; z_{p+1}  \, + \,  i c_{\tau_p }  } \; \cdot \; \f{\dd^n z}{ (2i\pi)^{n}} \;. 
\label{ecriture representation par integrale de boucle non modifiee}
\end{multline}
Focusing on the $z_p$-dependent part of the integration, we see that the part involving $I_N$ is holomorphic inside of 
$\Ga(J)$ (the singular part is due to $\chi^{-1}(z_p)$). Hence, by closing the loop we see that this part 
of the integral does not contribute to the value of the integral. In fact, one can replace $I_N$ by any
holomorphic function of $z_p$ inside of the loop $\Ga(J)$. Hence, we replace $I_N$ by $0$. 
Repeating this procedure for every $p$, we get 
\beq
I_{\{\tau_p\}}^{(n)} \; = \;  (-1)^n\Oint{\Ga(J)}{} \f{\dd^n z}{ (2i\pi)^{n}} 
\pl{p=1}{n}  \f{ \big( \bs{e}_{v_{\tau_{p-1}}},   \chi^{-1} (z_p)  \cdot 
\chi(z_{p+1} \, - i  c_{\tau_p} ) \cdot \bs{e}_{\tau_p} \big) } 
{ z_p \; - \; z_{p+1}  \, + \,  i c_{\tau_p }  } \;.
\enq

It only remains to identify the sums over $\tau_p \in \mc{E}_N$, $p=1,\dots, n$ as a matrix trace over the additional matrix indices
of the integral operator $\mc{M}$ so as to obtain the Fredholm determinant $\ln\det \big[I+\kappa \mc{M} \big]$. Taking the 
exponent, one obtains an equality, valid in some open neighborhood of $0$, between two entire functions. They are thus equal everywhere, 
and, in particular, at $\kappa=1$.   \qed

In fact, one can even derive an alternative representation for $\underset{\Ga(J)}{\det}\big[ I+\mc{M} \big]$ by deforming the loop $\Ga(J)$. 
The new matrix integral operator acts on the whole real axis.

\begin{prop}
One has the representation
\beq
\underset{\Ga(J)}{\det}[I+ \mc{M} ] \; = \; \underset{\R}{\det}[ I+ \mc{N} ] \;, 
\enq
where $\mc{N}$ is a $N\times N$ matrix integral operator on $L^2(\R)\otimes \Cx^N$  whose entries are given by 
\beq
\mc{N}_{k\ell}(\la,\mu) \; = \;   \e{sgn}(c_k) \cdot \ga_k \cdot 
 \f{ \big( \bs{e}_{v_k}\, ,  \, \big[ I_N \, - \, \chi^{-1}(\la+i\tf{c_k}{2}) \cdot \chi(\mu-i\tf{c_{\ell}}{2}) \big] 
 				\cdot \bs{e}_{\ell}  \big) }
 		{ 2i\pi \big( \la-\mu  + i \tf{(c_k+c_{\ell})}{2} \big) } \;. 
\label{ecriture matrice N}
\enq
\end{prop}

We do stress that the obtained operator is non singular in that, should $c_k+c_{\ell}=0$ for some $k,\ell$, then 
the simple zero of the denominator at $\la = \mu$ is compensated by a zero of the numerator. 

\Proof 

The starting point for obtaining the representation \eqref{ecriture matrice N} is
to recast \eqref{ecriture representation par integrale de boucle non modifiee} in a slightly different manner with the help of
contour deformations. Hence, starting from \eqref{ecriture representation par integrale de boucle non modifiee}, 
in the integral relative to $z_p$, we now replace $I_N$
by $\chi^{-1}(z_{p+1}-i c_{\tau_p})$ -instead of $0$-. 
Note that such a replacement also does not affect the analyticity properties of the integrand in $z_{p+1}$
belonging to the inside of the sufficiently small loop $\Ga(J)$. Repeating this procedure for every $p$, we get 
\beq
I_{\{\tau_p\}}^{(n)} \; = \; \Oint{\Ga(J)}{} 
\pl{p=1}{n}  \f{ \Big( \bs{e}_{v_{\tau_{p-1}}}, \Big[I_N -   \chi^{-1} (z_p)  \cdot 
\chi(z_{p+1} \, - i  c_{\tau_p} ) \Big]\cdot \bs{e}_{\tau_p} \Big) } 
{ 2i\pi (z_p \; - \; z_{p+1}  \, + \,  i c_{\tau_p } )  } \cdot \dd^n z \;.  
\enq

We now pass on to deforming the loop $\Ga(J)$. For this it is enough to consider the model integral
\beq
\mc{I}_{p} \; =   \;  \Oint{ \Ga(J) }{}  f_p(z_p)   \cdot   \dd z_p
\label{definition integrale modele}
\enq
where 
\beq
f_p(z_p) \; = \;  \f{ 
\Big( \bs{e}_{v_{\tau_{p-2}}}, \Big[I_N -   \chi^{-1} (z_{p-1})  \cdot \chi(z_{p} \, - i  c_{\tau_{p-1}} ) \Big]
				\cdot \bs{e}_{\tau_{p-1}} \Big)  
\cdot 
\Big( \bs{e}_{v_{\tau_{p-1}}}, \Big[I_N -   \chi^{-1} (z_p)  \cdot \chi(z_{p+1} \, - i  c_{\tau_p} ) \Big]\cdot \bs{e}_{\tau_p} \Big) } 
{( z_{p-1} \; - \; z_p  \, + \,  i c_{\tau_{p-1} } ) \cdot ( z_p \; - \; z_{p+1}  \, + \,  i c_{\tau_{p}} ) }
\enq
and then apply the results recursively. 
We deform the integration contour in \eqref{definition integrale modele} as
\beq
\Oint{ \Ga( J ) }{}  \quad \hookrightarrow \quad  \sul{ \veps = \pm }{} \; \Int{ \R_{-\veps}+i\veps 0^+ }{} \qquad 
\e{where} \qquad \R_{\veps} \equiv \intff{-\veps \infty}{ \veps \infty} \;. 
\enq
Let $s_p = \e{sgn}(c_{\tau_p})$. The integrand is analytic for $z_{p} \in \mathbb{H}_{-s_{p-1}}$ since the apparent pole
at $z_p = z_{p+1} - i c_{\tau_{p}}$ - which may or may not belong to $\mathbb{H}_{-s_{p-1}}$- is canceled by a zero of the numerator. 
Furthermore, $f_p(z_p) \; = \;  \e{O}(z_p^{-2})$ when $z_p \tend \infty$. 
Hence, by deforming the contour 
$\R_{s_{p-1}} - i s_{p-1}0^+ $ up to $\R_{s_{p-1}} - i s_{p-1} \times +\infty $, 
we can drop the corresponding contribution, leading to 
\beq
\mc{I}_{p} \; =   \hspace{-3mm}  \Int{ \R_{-s_{p-1} } + is_{p-1} 0^+  }{} \hspace{-3mm} f_p(z_p)  \cdot  \dd z_p  \;. 
\enq
 Repeating these manipulations in every integral arising in $I^{(n)}_{\{\tau_p\}}$, we get 
\beq
I_{\{\tau_p\}}^{(n)} \; =   \;  
\pl{p=1}{n} \bigg\{ \Int{ \substack{ \R_{-s_{p-1} }  \\ \;\; +is_{p-1} 0^+ }  }{}  \hspace{-3mm} \dd z_p  \bigg\}
\cdot \pl{p=1}{n}  \f{ \Big( \bs{e}_{v_{\tau_{p-1}}}, \Big[I_N -   \chi^{-1} (z_p)  \cdot 
\chi(z_{p+1} \, - i  c_{\tau_p} ) \Big]\cdot \bs{e}_{\tau_p} \Big) } 
{ 2i\pi (z_p \; - \; z_{p+1}  \, + \,  i c_{\tau_p } )  } \; . 
\enq
The integrand is holomorphic in the multidimensional strip in $\Cx^n$
\beq
\mc{S} \; = \; \Big\{ (z_1,\dots, z_n) \in \Cx^n \; : \; 0 <   s_{p-1}\Im(z_p) <  s_{p-1}\, c_{\tau_{p-1}}   \; \; , p=1,\dots, n   \Big\}
\enq
in that the numerator is holomorphic on this strip whereas the apparent simple poles due to the 
denominator are, in fact, canceled by the numerator's zeros. Thence,  
we can deform the integration contour in $z_p$
\beq
\e{from } \qquad \R_{-s_{p-1} } + is_{p-1} 0^+   \qquad \e{to}  \qquad 
\R_{-s_{p-1} } + i \tf{ c_{\tau_{p-1}} }{2}  \;. 
\enq
This leads to 
\beq
I_{\{\tau_p\}}^{(n)} \; =   \;  
\Int{ \R^n  }{}  
\pl{p=1}{n}\f{ \Big(\bs{e}_{v_{\tau_{p-1}}} \, ,  \, 
\big[ I_N \,  - \,     \chi^{-1}  (z_{p} \, +  \,  i\tf{c_{\tau_{p-1}}}{2} ) 
			\cdot  \chi(z_{p+1} \, -\,  i \tf{ c_{\tau_p} }{2} ) \big]  \bs{e}_{v_{\tau_p}} \Big)   } 
{ - \e{sgn}(c_{\tau_p}) \cdot \big( z_{p} \; - \; z_{p+1}  \, + \,  i \tf{( c_{\tau_{p-1}} + c_{\tau_{p}}) }{2} \big)   }   \;  \f{\dd^n z }{(2i\pi)^{n}} \;. 
\enq

It then only remains to recognize the sum over $\tau_p \in \mc{E}_N$, $p=1,\dots, N$ 
as the trace associated with the products over the matrix indices. The result then follows by an analytic continuation in $\kappa$,
just as in the previous proposition.  \qed

\section{Application to the two-shifted generalized sine kernel} 
\label{Section Application of Results}

In this section, we apply our general formalism to the case of the so-called two-shifted generalized sine kernel
whose kernel has been defined in \eqref{definition noyau type shift}. In such a way, we prove theorem \ref{Theorem principal}. 
The kernel \eqref{definition noyau type shift} is of the form \eqref{definition noyau avec shifts} 
with $N=2$, $v_a=a$, $\ga_1=\ga_2=1$ and $c_1=-c$, $c_2=c$ and the two-dimensional vectors 
$\bs{E}_{R}(\la)$ and $\bs{E}_{L}(\la)$ being expressed as: 
\beq
 \bs{E}_{L}^{\bs{T}}(\la) \; = \;   \f{ F(\la) }{ 2i\pi }  \big( - e^{-1}\!(\la) \; , \;   e(\la) \big)  \qquad \e{and} \qquad 
 \bs{E}_{R}^{\bs{T}}(\la) \; =  \;   \big( e(\la) \; , \; e^{-1}(\la) \big)  \;. 
\enq

It has been established in \cite{KozKitMailSlaTerRHPapproachtoSuperSineKernel} that, when $x\tend +\infty$,
given any open and relatively compact neighborhood $O$ of $\intff{a}{b}$, 
the unique solution $\chi$ to the Riemann--Hilbert problem associated with the 
operator $\wt{S}$ \eqref{definition operateur S tilde noyau}  takes, on  $\Cx \setminus \ov{O}$, the form 
\beq
\chi(\la)  \; =  \; \Pi(\la) \cdot \a^{-\sg_3}(\la) \qquad \e{with} \qquad 
\norm{ \Pi -I_2}_{L^{\infty}\big( \Cx \setminus \ov{O} \big) } \; = \; \e{O}\Big( \f{1}{x} \Big) \;. 
\enq
Here, $\Pi$ is some holomorphic matrix in $\Cx \setminus \ov{O}$ and the function $\a$
is as defined by \eqref{definition fonction alpha}.

Thus we choose some relatively compact neighborhood $O$ of $\intff{a}{b}$ that is small enough 
and take $\Ga$  to be a small counterclockwise loop around $\intff{a}{b}$, lying entirely in 
$\big\{\Cx \setminus \ov{O} \big\} \cap \big\{ | \Im(z) | < \tf{c}{2} \big\}$. 
Then, we introduce two integral operators $I\; + \; \mc{M}^{(0)}$ and $I\; + \; \mc{M}$ 
on $L^2\big( \Ga \big) \otimes \Cx^2$ having kernels 
\beq
\mc{M}(\la,\mu) \; = \;  \pa{  \ba{cc}  
  \f{ \big(  \bs{e}_1 \, ,  \chi^{-1}(\la) \cdot \chi(\mu + ic) \cdot  \bs{e}_1 \big) }{ 2i\pi ( \la-\mu -i c ) }  & 
   \f{ \big(  \bs{e}_1 \, , \,  \chi^{-1}(\la) \cdot \chi(\mu - i c) \cdot \bs{e}_2 \big)  }{ 2i\pi( \la-\mu +i c) } \vspace{2mm} \\
 \f{  \big(  \bs{e}_2 \, , \,  \chi^{-1}(\la ) \cdot \chi(\mu+ic) \cdot \bs{e}_2 \big)  }{ 2i\pi (\la-\mu -i c) }   &
  \f{ \big(  \bs{e}_2 \, , \, \chi^{-1}(\la) \cdot \chi( \mu - ic) \cdot  \bs{e}_2 \big)  }{  2i\pi(\la-\mu  + i c) }      \ea } \;. 
\label{ecriture matrice M pour GSK}
\enq
and 
\beq
\mc{M}^{(0)}(\la,\mu) \; = \;  \pa{  \ba{cc}  
  \f{ \a(\la) \cdot \a^{-1}(\mu+ic) }{ 2i\pi \cdot ( \la-\mu -i c ) }  &    0   \\
0   &
 \f{ \a^{-1}(\la) \cdot \a(\mu-ic) }{  2i\pi \cdot ( \la-\mu  + i c ) }      \ea } \;. 
\label{ecriture approximant M0 pour GSK}
\enq

$\mc{M}(\la,\mu)$ and $\mc{M}^{(0)}(\la,\mu)$ being bounded on $\Ga^2$,  and $\Ga$ being compact, 
the integral operators $\mc{M}$ and $\mc{M}^{(0)}$ are trace class. 
It thus follows from standard estimates for Fredholm determinants (\textit{cf} \cite{SimonsInfiniteDimensionalDeterminants})
that, for some universal constant $C$, 
\beq
\Big| \det\big[ I+\mc{M}^{(0)} \big] \; - \; \det \big[  I+\mc{M} \big]   \Big|  \; \leq \; C
\norm{ \mc{M}^{(0)} - \mc{M} }_1 \;, 
\enq
with $\norm{ \cdot }_1$ being the trace class norm. Hence, the uniform bounds for $\Pi-I_2$ on $\Ga$ ensure that 
\beq
\det\big[ I+\mc{M} \big] \; = \;  \det\big[ I+\mc{M}^{(0)} \big] \cdot   \bigg\{ 1 \; + \; \e{O} \Big( \f{1}{x} \Big)  \bigg\} \; .  
\enq

It is then enough to apply the propositions 
\ref{Proposition factorization determiant V en terme RHP sovable simple det V tile et det W}
and \ref{Proposition factorisation Det W en terme RHP data det M} specialized to the setting
associated with the kernel $S$ \eqref{definition noyau type shift} and 
 choose the loop arising in proposition \ref{Proposition factorisation Det W en terme RHP data det M} to
coincide with $\Ga$ (what is indeed possible). \qed

\section*{Conclusion}

In this paper, we have developed an effective technique for dealing with certain classes of integrable integral 
operators with kernels of shift-type. This allowed us to compute, in the case of the simple example 
of the generalized sine kernel issued two-shift integrable integral operator, the large-$x$ asymptotic behavior of 
its Fredholm determinant. The main achievement of this work is a technique that allows one to circumvent,
in a very simple way, the handling of operator valued Riemann--Hilbert problems. 
It is still interesting to compare the approach we have developed in this
paper with the operator valued Riemann-Hilbert technique of \cite{ItsSlavnovNLSTimeAndSpaceCorrDualFields}. We are going to
address this issue in a forthcoming publication.

\section*{Acknowledgements}

K. K. K. is supported by CNRS. He also acknowledges a funding from 
the grant PEPS-PTI "Asymptotique d'int\'{e}grales multiples". 
The work of A. R. Its was supported in
part by NSF grant DMS-1001777.
K. K. K. would like to thank the mathematics department of IUPUI for its warm hospitality and financial support during 
his visit there when this work has been carried out.

%%%%%%%%%%%%%%%%%%%%%%%%%%%%%%%%%%%%%%%%%%%%%%%%%%%%%%%%%%%%%%%%%%%%%%%%%%%%%%%%%%%%%%%%%%%%%%%%%%%%%%%%%%%%%%%%%%%%%%%%%%%%%%%%%%%%%%%%%%%%%%%%%%%%%%%%%%%%%%%%%%%%%%%%%%%%%%%%%%%%%%%
%%%%%%%%%%%%%%%%%%%%%%%%%%%%%%%%%%%%%%%%%%%%%%%%%%%%%%%%%%%%%%%%%%%%%%%%%%%%%%%%%%%%%%%%%%%%%%%%%%%%%%%%%%%%%%%%%%%%%%%%%%%%%%%%%%%%%%%%%%%%%%%%%%%%%%%%%%%%%%%%%%%%%%%%%%%%%%%%%%%%%%%

\end{document}